\documentclass[a4paper]{amsart}
\usepackage{amssymb}
\usepackage{amsfonts}
\usepackage{graphicx}
\usepackage{amsmath}
\usepackage{hyperref}
\usepackage{url}

\newtheorem{theorem}{Theorem}

\newtheorem{proposition}[theorem]{Proposition}

\begin{document}
\title{The Structure of Z-related Sets}
\author{Franck Jedrzejewski, Tom Johnson}
\maketitle

\begin{abstract}
The paper presents some new results on Z-related sets obtained by
computational methods. We give a complete enumeration of all Z-related sets
in $\mathbb{Z}_{N}$ for small $N$. Furthermore, we establish that there is a
reasonable permutation group action representing the Z-relation.

\medskip

\end{abstract}



\bigskip

The concept of Z-relation has been discussed in a systematic way by Allen
Forte \cite{For1977}. In the 1940s, the concept appears in crystallography
under the name of \textquotedblleft homometric point sets\textquotedblright\
in the paper of A.L.\ Patterson \cite{Pat1944}. In these years, homometric
sets have been extensively studied by mathematicians and crystallographers.\
R.K.\ Bullough pointed out some general theorems \cite{Bul1961} until J.\
Rosenblatt rephrased and extended the problem in a new algebraic framework.\
In music theory, the problem was studied by S.\ Soderberg in 1995 \cite%
{Sod1995}, and by J.\ Goyette in 2012 \cite{Goy2012}.\ J.\ Mandereau \emph{%
et al.}\ \cite{Man2011a}, \cite{Man2011b} gave a new starting point in\ 2011.

\section{Z-relation and homometry}

A pitch class set is a set of pitches where octaves are equivalent, and
enharmonically equivalent pitches are identified.\ In music theory, the
Z-relation links two different set classes with the same interval vector.\ $%
N $-tone equal temperament is represented by pitch classes $\mathbb{Z}%
_{N}=\{0,1,...,N-1\}$.\ For a given set $A\subset \mathbb{Z}_{N},$ the \emph{%
interval vector} measures the number of ways the interval $n$ can be spanned
between members of $A$.%
\begin{equation*}
\mathrm{iv}(A)(n)=\mathrm{ifunc}(A,A)(n)=\#\{(a,b)\in A\times A,~b-a=n\}
\end{equation*}%
The \emph{interval function} of two sets $A$ and $B$, introduced by David
Lewin in \cite{Lew1987}, is the number of times any $k$ in $A$ has its $n$%
-transpose in $B$, namely:%
\begin{equation*}
\mathrm{ifunc}(A,B)(n)=\underset{k}{\sum }\mathbf{1}_{A}(k)\mathbf{1}%
_{B}(n+k)
\end{equation*}%
where $\mathbf{1}_{A}$\ is the indicator function defined by $\mathbf{1}%
_{A}(k)=1$ if $k\in A$ and $\mathbf{1}_{A}(k)=0,$ otherwise. If the sets $A$
and $B$ are equal, the interval function is the interval vector. The \emph{%
interval content} is the set of first digits of the interval vector, ($%
\mathrm{iv}(A)(1),,..,\mathrm{iv}(A)([N/2]))$, except that the last one is
divided by 2 if $N$ is even. For example, the interval vector of the set
\{0, 1, 3, 4, 7, 9\} is (6, 2, 2, 4, 3, 2, 4, 2, 3, 4, 2, 2), and its
interval content is (2,2,4,3,2,2).

Two sets $A$ and $B$ of $\mathbb{Z}_{N}$ are said to be $Z$-related if they
have the same interval content:%
\begin{equation*}
(A~\mathcal{Z}_{N}~B)\Leftrightarrow \mathrm{ic}(A)=\mathrm{ic}(B)
\end{equation*}%
In other words, $A$ and $B$ share the same interval function, or the same
interval vector. Since transposing or inverting do not change the interval
content, we get a lot of trivially $Z$-related sets.\ To avoid trivial
cases, we consider set classes up to transposition and inversion, that is
under the action of the dihedral group.\ The first well-known example was
given by Lino Patterson in 1944:

\begin{equation*}
\{0,3,4,5\}~\mathcal{Z}_{8}~\{0,4,5,7\}
\end{equation*}%
If we draw the two sets of points on a circle representing $\mathbb{Z}_{8}$,
we will see that the distances between the points are the same but the
points are in different arrangements.\ The two sets share the same interval
content $ic=2121$ (meaning from left to right that there are 2 pairs with
distance 1, 1 pair with distance 2, 2 pairs with distance 3, and 1 pair with
distance 4). In dealing with finite abelian groups $\mathbb{Z}_{N}$, it will
be convenient to use polynomial notation. The polynomial representing the
set $A$ of $\mathbb{Z}_{N}$\ is the polynomial whose exponents are the
elements of $A$:%
\begin{equation*}
A(x)=\underset{a\in A}{\sum }x^{a}
\end{equation*}%
The \textit{reflection} of $A$ is the polynomial whose exponents are the
inversion $I(A)$ of $A$ mod $N$%
\begin{equation*}
A^{\ast }(x)=\underset{b\in I(A)}{\sum }x^{b}
\end{equation*}%
$I$ denotes the inversion $I(x)=-x~\mathop{\rm mod}\nolimits~N$. By
definition, the \textit{Patterson function} of the set $A$ is the
autocorrelation function given by the convolution product%
\begin{equation*}
F(x)=A\ast A^{\ast }(x)=A(x)A(x^{-1})
\end{equation*}%
modulo $(x^{N}-1)$.\ The Patterson function is an equivalent of the interval
function, since we have the following result.\ 

\begin{proposition}
The coefficients of the Patterson function are equal to the components of
the interval function.%
\begin{equation*}
A\ast A^{\ast }(x)=\underset{c_{i}\in \mathrm{ifunc}(A)}{\sum }c_{i}x^{i}
\end{equation*}
\end{proposition}

Two sets $A$ and $B$\ of $\mathbb{Z}_{N}$\ are said to be \textit{homometric}
if they have the same Patterson function. The musical concept of Z-related
sets coincides with the crystallographic notion of homometric sets.

\bigskip

\noindent \textit{Example}\emph{.}\ For $N=12$ and $A=\{0,2,3,5\},$ it is
easy to verify the previous result.{\footnotesize 
\begin{eqnarray*}
A(x) &=&1+x^{2}+x^{3}+x^{5},\qquad A^{\ast }(x)=1+x^{10}+x^{9}+x^{7} \\
F(x) &=&A(x)A(x^{-1})=4+x+2x^{2}+2x^{3}+x^{5}+x^{7}+2x^{9}+2x^{10}+x^{11} \\
\mathrm{ifunc}(A) &=&(4,1,2,2,0,1,0,1,0,2,2,1)
\end{eqnarray*}%
}

\section{Some general theorems}

In $\mathbb{Z}_{N}$, the set $A$ is trivially homometric under the action of
the dihedral group, that is under all inversions $I_{n}(A)$ and
transpositions $T_{n}(A)$, where inversions are defined by $I_{n}(x)=-x+n~%
\mathop{\rm mod}\nolimits~N$, and transpositions are defined by $%
T_{n}(x)=x+n~\mathop{\rm mod}\nolimits~N.$

Moreover, let $A$ and $B$ be subsets of $\mathbb{Z}_{N}$, then $A$ and $B$
are said to be \emph{trivially homometric} if they belong to the same orbit
under the action of the dihedral. If $A$ and $B$ are homometric but not
trivially homometric, then they are \emph{strictly homometric}.\ Thus the
Z-relation is equivalent to strict homometry.

In 1944, Patterson \cite{Pat1944} established the following two results:

\noindent (1) If two subsets of a regular $N$-gon are homometric then their
complements are. 
\begin{equation*}
A~\mathcal{Z}_{N}~B\Leftrightarrow A^{c}~\mathcal{Z}_{N}~B^{c}
\end{equation*}%
(2) Every $N$-point subset of a regular $2N$-gon is homometric to its
complement.

In terms of music theory, the $N$-note set $A$ of the $2N$-tone equal
temperament is in Z-relation with its complement $A^{c}$. For example, for $%
N=4$, the 4-note set $A=\{0,1,3,5\}$ of $\mathbb{Z}_{8}$ and its complement $%
A^{c}=\{2,4,6,7\}$ are in Z-relation since

\begin{equation*}
\mathrm{ic}(A)=(1,2,2,1)=\mathrm{ic}(A^{c})
\end{equation*}

If we consider now the multiplication by an integer $m$ modulo $N$, it is
easy to show that the Z-relation is stable by multiplication.\ More
precisely, let $A$\ be a set in $\mathbb{Z}_{N},$ $m$ an integer such that $%
\mathrm{gcd}(m,N)=1$ and $m\neq 1,N-1.$ Then the interval content $\mathrm{ic%
}(A)$ is in general different from $\mathrm{ic}(M_{m}(A))$, where $M_{m}$
denotes the multiplication by $m$ modulo $N,$ $M_{m}(x)=mx~\mathop{\rm mod}%
\nolimits~N$. We have%
\begin{equation*}
A~\mathcal{Z}_{N}~B\Longrightarrow M_{m}A~\mathcal{Z}_{N}~M_{m}B
\end{equation*}

\bigskip

\noindent \textit{Example}. In the usual temperament ($N=12$), since the two
complementary sets $A=\{0,1,2,3,5,6\}$ and $B=A^{c}=\{0,1,2,3,4,7\}$ are in
Z-relation with interval content $\mathrm{ic}(A)=433222,$ we get a new
homometric pairs by mutliplying each set by 5.\ Thus, $M_{5}A=\{0,1,3,5,6,10%
\}$ and $M_{5}B=\{0,3,5,8,10,11\}$ are Z-related with a new interval content
equal to $\mathrm{ic}(M_{5}A)=233242.$

\bigskip

\noindent \textit{Remark}. Unfortunately, not all Z-related pairs are coming
from complement nor multiplication, $A~\mathcal{Z}_{N}~A^{c}$ or $A~\mathcal{%
Z}_{N}~M_{m}A$. For example, for $N=18$, the sets $A=\{0,1,2,3,4,7,8,14,16\}$
and $B=\{0,1,2,3,5,6,7,9,13\}$ are Z-related but $B$ is neither the
complementary set of $A$, nor a multiple set of $A$ ($B\neq M_{k}A$\ for $%
k=5,7,11,13$).

\bigskip

By adding the transposed set to the initial set, we get some new theorems:

\begin{theorem}
Let A and B be two Z-related sets in $\mathbb{Z}_{N}$, and $T_{N}$ the
transposition in $\mathbb{Z}_{2N}$ ($T_{N}(x)=x+N$ $\mathop{\rm mod}\nolimits
$ $2N$), then we have%
\begin{equation*}
A~\mathcal{Z}_{N}~B\Longrightarrow (A\cup T_{N}A)~\mathcal{Z}_{2N}~(B\cup
T_{N}B)
\end{equation*}
\end{theorem}

Or more generally,

\begin{theorem}
Let A and B be two Z-related sets in $\mathbb{Z}_{N}$, and $T_{j}$ the
transpositions in $\mathbb{Z}_{Nm}$, then we have%
\begin{equation*}
A~\mathcal{Z}_{N}~B\Longrightarrow (A\cup T_{N}A\cup ...\cup T_{N(m-1)}A)~%
\mathcal{Z}_{Nm}~(B\cup T_{N}B\cup ...\cup T_{N(m-1)}B)
\end{equation*}
\end{theorem}

The same can be done with the multiplication.\ 

\begin{theorem}
Let $A~\mathcal{Z}_{N}~B$ and $m$ be an integer such that $\mathrm{gcd}%
(N,m)=1$ and $m\neq 1,N-1$, $M_{m}$ be the multiplication and $T_{j}$ the
transpositions in $\mathbb{Z}_{Nm},$ then we have%
\begin{equation*}
(M_{m}A\cup T_{1}M_{m}A\cup ...\cup T_{m-1}M_{m}A)~\mathcal{Z}%
_{Nm}~(M_{m}B\cup T_{1}M_{m}B\cup ...\cup T_{m-1}M_{m}B)
\end{equation*}
\end{theorem}

\bigskip

\noindent \textit{Example}. Starting from the homometric pair $\{0,1,3,4\}~%
\mathcal{Z}_{8}~\{0,1,2,5\},$ the multiple relation of $M_{3}A=\{0,3,9,12\}$
and $M_{3}B=\{0,3,6,15\}$ leads to the new pair:

\begin{center}
$\{0,1,2,3,4,5,9,10,11,12,13,14\}~\mathcal{Z}_{24}~%
\{0,1,2,3,4,5,6,7,8,15,16,17\}$
\end{center}

\bigskip

In 2008, O'Rourke, Taslakian and Toussaint \cite{Oro2008} gave a new
procedure for generating homometric pairs called the \emph{Pumping Lemma,}
based on adding points around isospectral vertices (see their paper for
details). We propose another new procedure starting with not one, but two
homometric pairs. In the general case, this procedure is lengthy and will be
published elsewhere. We consider here an example.\ Choose two suitable
pairs.\ For example $\{0,1,3,4,6\}~\mathcal{Z}_{10}~\{0,1,2,4,7\}$ and $%
\{0,1,4,5,7\}~\mathcal{Z}_{12}~\{0,1,2,5,8\}.\ $Place the first of each pair
on a circle of 10 points (Fig. 1, left), and add the second of each pair on
the circle as shown on Fig. 1 (center).\ On each of the two points connected
by a diameter, we can add $k$ points, leading to the new relation (Fig. 1
right)%
\begin{equation*}
\{0,1,3+k,4+k,6+k\}~\mathcal{Z}_{10+2k}~\{0,1,2,4+k,7+k\}
\end{equation*}%
Remark that the number of added points is always even, even if they are
several diameters.

\bigskip

\begin{figure}[!h]
\centering
\includegraphics[height=4cm]{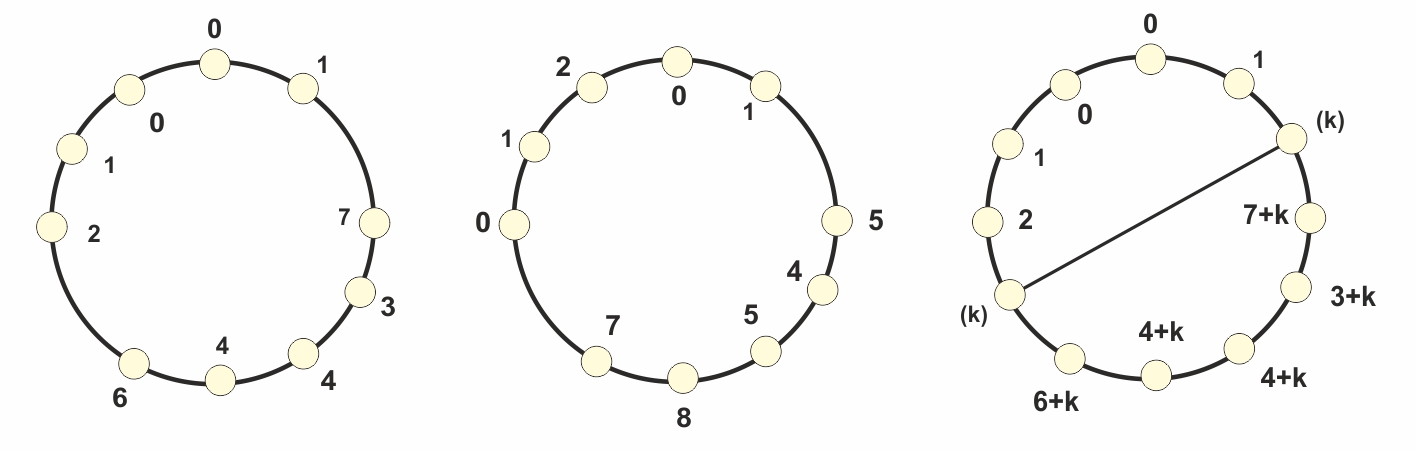}
\caption{The two interlaced pairs}
\label{fig:ex1}
\end{figure}

\bigskip

The aim of the pumping lemma or of this procedure is to find a way of
establishing general relations, such as this one, established empirically.\
For $n\geq 5$, we have%
\begin{eqnarray*}
&&\{0,1,n-2,n-1,n+1\}~\mathcal{Z}_{2n}~\{0,1,2,n-1,n+2\} \\
&&\{0,1,2,n-2,n+1\}~\mathcal{Z}_{2n}~\{0,1,3,n-1,n\}
\end{eqnarray*}%
Rosenblatt gave a complete classification of Z-related pairs of length 4:

\bigskip

If $A~\mathcal{Z}~B$ with $\mathrm{card}(A)=\mathrm{card}(B)=4$ then $A$ and 
$B$ are of the following two types:

(i) In $\mathbb{Z}_{4n}$, $\exists a\in \{1,2,...,n-1\},~n\geq 2,$

\begin{equation*}
A=\{0,a,a+n,2n\}~\mathcal{Z}_{4n}~B=\{0,a,n,2n+a\}
\end{equation*}

(ii) In $\mathbb{Z}_{13n}$,

\begin{equation*}
A=\{0,n,4n,6n\}~\mathcal{Z}_{13n}~B=\{0,2n,3n,7n\}
\end{equation*}%
But the question of finding the classification of all Z-related pairs for $%
\mathrm{card}(A)\geq 5$ is still open (see \cite{Alt2000, Cal2008}).

\section{Enumeration of Z-related sets}

Enumeration of Z-related sets was studied by Patterson, and later by others,
including B\"{u}rger \cite{Bue1977} and Chieh \cite{Chi1979}.\ Patterson 
\cite{Pat1944}\ was the first to demonstrate the existence of $t$-uples of
homometric sets. For $N\leq 12,$ homometric sets are well known.\ Some
progress has been made for \textquotedblleft
complementary\textquotedblright\ homometric pairs with $N$ even. But from $%
N\geq 16,$ non-complementary homometric pairs appear.\ In \cite{Lew1981}
David Lewin rediscovered, and brought to the attention of the
music-theoretic community, the existence of the Z-triples for $N=16$, and
Jon Wild \cite{Wil1996} has previously entabulated all homometric tuplets
for $N\geq 31$, with tuplets of 15-note sets with as many as 16 members in $%
\mathbb{Z}_{30}$.\ But enumeration is still an open problem. Our computation
led to the following table.\ For each $k$-note sets, the table gives the
number of (distinct) interval vectors in $\mathbb{Z}_{N}$ for which
non-trivial homometric tuples exist (sets are considered up to inversion and
transpositions). For example, for $N=16$, there are 31 non-trivial
homometric tuples of $k=6$ notes. (In fact, 28 homometric pairs +\ 3\
homometric triplets).

\bigskip

\begin{equation*}
\begin{tabular}{c|cccccccccc}
\hline
$k\backslash N$ & $8$ & $10$ & $12$ & $13$ & $14$ & $15$ & $16$ & $17$ & $18$
& $19$ \\ \hline
4 & 1 & 0 & 1 & 1 & 0 & 0 & 2 & 0 & 0 & 0 \\ 
5 & -- & 3 & 3 & 0 & 6 & 5 & 10 & 0 & \emph{14} & 0 \\ 
6 & -- & -- & 15 & 2 & 6 & 25 & \emph{31} & 16 & \emph{62} & 21 \\ 
7 & -- & -- & -- & -- & 48 & 10 & \emph{44} & 24 & \emph{134} & 57 \\ 
8 & -- & -- & -- & -- & -- & -- & \emph{180} & 52 & \emph{150} & 90 \\ 
9 & -- & -- & -- & -- & -- & -- & -- & -- & \emph{572} & 156 \\ \hline
\end{tabular}%
\end{equation*}

\bigskip

Italics indicate the existence of tuples ($t>2$). For $N=12$, there are 19
homometric pairs (1+3+15).

\bigskip

\begin{equation*}
\begin{tabular}{c|ccccc}
\hline
$k\backslash N$ & $20$ & $21$ & $22$ & $23$ & $24$ \\ \hline
4 & 2 & 0 & 0 & 0 & 3 \\ 
5 & 22 & 0 & 20 & 0 & \emph{31} \\ 
6 & \emph{98} & 96 & \emph{60} & 33 & \emph{275} \\ 
7 & \emph{191} & 220 & \emph{335} & 110 & \emph{676} \\ 
8 & \emph{535} & 282 & \emph{575} & 429 & \emph{2532} \\ 
9 & \emph{565} & \emph{1062} & \emph{1425} & \emph{814} & \emph{5112} \\ 
10 & \emph{2106} & 613 & \emph{1550} & 1144 & \emph{7715} \\ 
11 & -- & -- & \emph{7390} & \emph{1375} & ? \\ 
12 & -- & -- & -- & -- & ? \\ \hline
\end{tabular}%
\end{equation*}

\bigskip

The first triple appears for $N=16$ and length 6:%
\begin{equation*}
\begin{tabular}{lll}
$\{0,1,2,4,6,9\},$ & $\{0,1,2,4,9,14\},$ & $\{0,1,3,5,7,8\}$%
\end{tabular}%
\end{equation*}%
and the first quadruple is obtained for $N=18$ and length 9 (in fact, there
are 54 quadruples). There are three quintuples for $N=24$, $k=10$, one
sextuple for $N=24$, $k=7$ and one octuple for $N=24$\emph{,} $k=9$: a very
impressive property of the quarter-tone universe.%
\begin{equation*}
\begin{tabular}{lll}
$\{0,1,2,4,6,9,12,16,17\}$ & $\{0,1,2,4,6,9,14,17,18\}$ & $%
\{0,1,2,4,8,9,12,14,17\}$ \\ 
$\{0,1,2,4,9,10,14,17,22\}$ & $\{0,1,2,4,9,14,16,17,20\}$ & $%
\{0,1,2,6,9,10,12,14,17\}$ \\ 
$\{0,1,3,5,7,8,13,16,17\}$ & $\{0,1,3,5,8,9,13,15,16\}$ & 
\end{tabular}%
\end{equation*}

\section{Group Action for Z-relation}

Since our interest is to give a fast algorithm of finding all homometric
sets, one way is to understand how they are organized and to answer the
question: Is there a non-trivial group action representing the Z-relation
(such that the orbits are the equivalence classes of the Z-relation)? If you
look at a subgroup of the linear group, you certainly do not find a
solution, as shown by Mandereau \emph{et al.}\ \cite{Man2011b}.\ But if you
look at a subgroup of the permutation group $\mathfrak{S}_{N}$, you will
find a nice non-trivial solution.\ Considering the set $\mathcal{B}$ of all
transpositions and inversions for a given $N$ of all homometric sets of
length $k$, the automorphism group of $\mathcal{B}$ is clearly a solution of
the problem. In other words, the group is a subgroup of $\mathfrak{S}_{N}$\
whose action on the power set of $\mathbb{Z}_{N}$ stabilizes $\mathcal{B}$.
To compute this automorphism group, the idea is to use a well-known object
in combinatorial block design known as the \emph{Levi graph }$\Gamma $.
Denoting the homometric sets $\mathcal{B=\{}B_{1},...,B_{u}\}$, the vertex
set of the Levi graph is $V=\{0,1,2,...,N-1,B_{1},...,B_{u}\}.$ The $N$
first elements of $V$ have one color and the $\mathcal{B}$\ elements have a
second color.\ The edge set is defined by: 
\begin{equation*}
\{(i,B_{j}):i\in B_{j},\text{ }i=0,1,...,N-1\}
\end{equation*}%
As is well known, the automorphism group of the Levi graph is the same as
the automorphism group of $\mathcal{B}.$ Thus, the problem is to compute the
automorphism group of a vertex-colored graph, which is done by a C program.
The information is provided in the form of a set of generators, the size of
the group and the orbits of the group.\ After this computation, the
difficulty is to give a nice representation of the orbits of the group.
Let's have a look at some examples.

\bigskip

\noindent \textit{Example 1}. The simplest example is to compute the
automorphism group of the homometric sets for $N=8$.\ As we know, there is
only one Z-related pair, but 16 transpositions and inversions. The pair

\begin{equation*}
\{0,1,2,5\}~\mathcal{Z}_{8}~\{0,1,3,4\}
\end{equation*}%
leads to a set $\mathcal{B}$ of 16 elements.\ The automorphism group is
generated by four generators {\small $a=(1,3)(2,6)(5,7),$ $b=(1,5)(3,7),$ $%
c=(0,1)(2,7)(3,6)(4,5),$ $d=(2,6)(3,7)$}.\ The group has only one orbit, and
is represented on Fig\ 2.\ The inner circle (dihedral class of the
homometric set \{0,1,2,5\}) is linked by permutation $d$ to the outer circle
(dihedral class of the homometric set \{0,1,3,4\}).

\bigskip

\begin{figure}[!h]
\centering
\includegraphics[height=7cm]{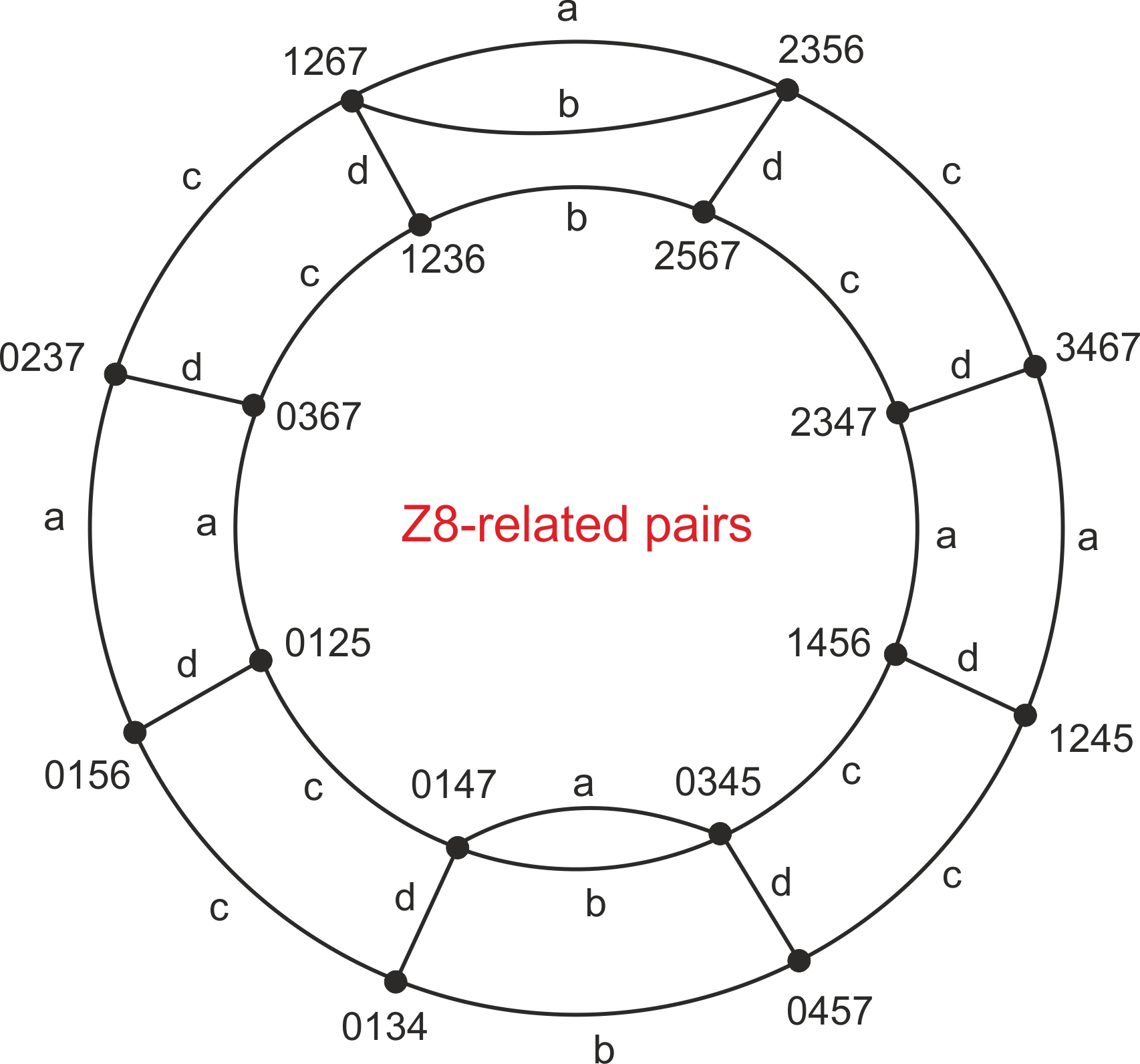}
\caption{Group of Z$_{8}$-related pairs}
\label{fig:ex2}
\end{figure}

\bigskip

\noindent \textit{Example 2}. For $N=12$, the automorphism group of the 48\
homometric sets generated by the pair $\{0,1,3,7\}$, $\{0,1,4,6\}$ of length
4, has six generators.

\begin{eqnarray*}
a &=&(3,9),~b=(4,10),~c=(5,11) \\
d &=&(2,5)(8,11),~e=(1,2)(4,5)(7,8)(10,11) \\
f &=&(0,1)(3,4)(6,7)(9,10)
\end{eqnarray*}%
Its representation is a rather complicated graph. The automorphism group of
the 108\ homometric sets of length 5 has three generators: 
\begin{eqnarray*}
a &=&(1,5)(2,10)(4,8)(7,11) \\
b &=&(1,7)(3,9)(5,11) \\
c &=&(0,1)(2,11)(3,10)(4,9)(5,8)(6,7)
\end{eqnarray*}%
Tom Johnson \cite{Tom2013}\ demonstrated that the homometric sets of length
5 can be graphed with less than six transformations, as one can see in the
concentric circles of Fig. 3.

\bigskip

\begin{figure}[!h]
\centering
\includegraphics[height=4.7in]{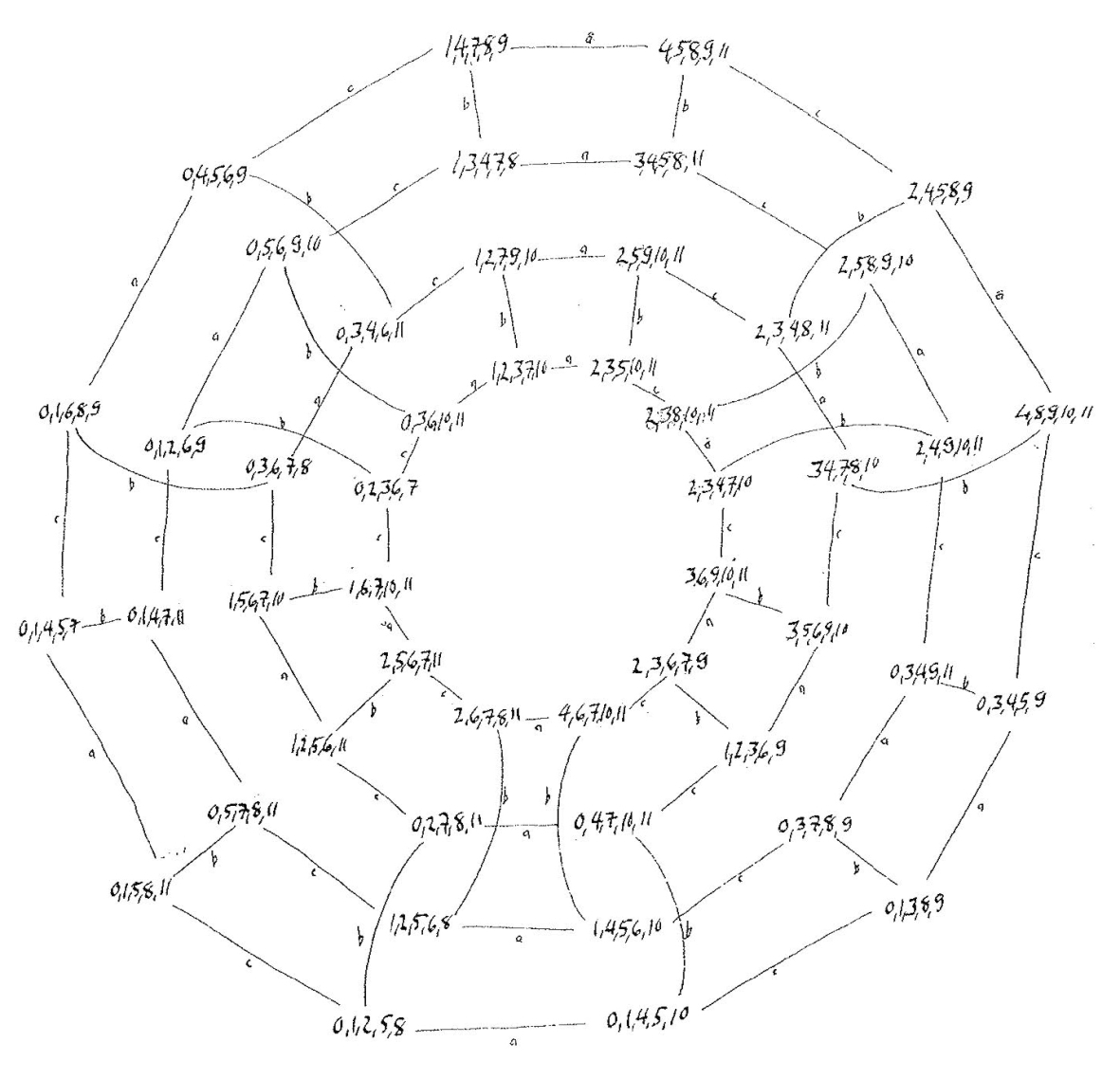}
\caption{Tom Johnson's drawing (Orbit of the Group of Z$_{12}$-related pairs
of length 5)}
\label{fig:ex3}
\end{figure}

\bigskip

The structure of this graph appears 6\ times in the 19 orbits of the
automorphism group of the 552 homometric pairs of length 6 whose generators
are

\begin{eqnarray*}
a &=&(2,10)(3,11)(4,8)(5,9) \\
b &=&(1,3)(2,10)(4,8)(5,11)(7,9) \\
c &=&(0,1)(2,3)(4,5)(6,7)(8,9)(10,11)
\end{eqnarray*}%
Tom Johnson has experimented with new transformations in a more musical way.
Starting with the 48 homometric four-note chords, he computed the
two-dimensional network of Table 1, described thus: For horizontal
transformations the note contained both in the M2 (Major second) and in the
m3 (minor third) moves a tritone. For odd vertical transformations the two
notes of the m3 move one place away from the m2. For even vertical
transformations, the notes of the M3 approach one another, becoming a M2, or
the two notes of the M2 separate into a M3.

\bigskip

\begin{center}
15ab\qquad 57ab\qquad 457b\qquad 145b

059b\qquad 569b\qquad 356b\qquad 035b

068b\qquad 0568\qquad 0256\qquad 026b

067a\qquad 0467\qquad 0146\qquad 016a

1679\qquad 1367\qquad 0137\qquad 0179

1578\qquad 1257\qquad 127b\qquad 178b

2478\qquad 1248\qquad 128a\qquad 278a

2368\qquad 0238\qquad 0289\qquad 2689

2359\qquad 239b\qquad 389b\qquad 3589

1349\qquad 139a\qquad 379a\qquad 3479

034a\qquad 049a\qquad 469a\qquad 346a

24ab\qquad 48ab\qquad 458a\qquad 245a

\medskip

Table 1\ -- The complete group of 48 homometric four-note chords
\end{center}

\bigskip

By considering transformations in musical terms, Tom Johnson opens new ways
of understanding homometric relationships and defines an automorphism group
at the same time.

\section{Conclusion}

The aim of this paper was to consider homometric sets and to compute the
number of interval vectors for each value of $N$. We establish some new
theorems helping us to compute Z-related sets.\ Moreover, we show that there
is a reasonable group representing all Z-related sets of a given length, for
a given N and compute this group explicitly for $N$ less than 12. However,
the question of the enumeration of homometric sets remains open.

\section*{acknowledgement}
Submitted to MCM\ 2013.\ The final publication will be publish by Springer
Lectures Notes in\ Computer Sciences and available at link.springer.com.
\url{http://link.springer.com/bookseries/558}.

\end{document}